# ON CONSTRUCTION OF THE SMALLEST ONE-SIDED CONFIDENCE INTERVAL FOR THE DIFFERENCE OF TWO PROPORTIONS

By Weizhen Wang[1]

For any class of one-sided $1-\alpha$ confidence intervals with a certain monotonicity ordering on the random confidence limit, the smallest interval, in the sense of the set inclusion for the difference of two proportions of two independent binomial random variables, is constructed based on a direct analysis of coverage probability function. A special ordering on the confidence limit is developed and the corresponding smallest confidence interval is derived. This interval is then applied to identify the minimum effective dose (MED) for binary data in dose-response studies, and a multiple test procedure that controls the familywise error rate at level $\alpha$ is obtained. A generalization of constructing the smallest one-sided confidence interval to other discrete sample spaces is discussed in the presence of nuisance parameters.

**1. Introduction.** We first focus on an important case, the binomial distribution. Let $X$ be a binomial random variable with $n$ trials and a probability of success $p_1$, denoted by $Bin(n, p_1)$, and let $Y$ be another independent $Bin(m, p_0)$. Their probability mass functions and cumulative distribution functions are denoted by $p_X(x; p_1, n)$, $p_Y(y; p_0, m)$, $F_X(x; p_1, n)$ and $F_Y(y; p_0, m)$, respectively. The goal of this paper is to construct the optimal one-sided $1-\alpha$ confidence interval of form $[L(X,Y), 1]$ for $p_1 - p_0$ and to discuss its application and a generalization to other discrete sample spaces. This type of interval is important when one needs to establish that $p_1$ is larger than $p_0$ by a certain amount. An immediate application is in a clinical trial where the goal is to determine if a treatment is "better" than a control with binary responses.

There are two general criteria used to evaluate the performance of a confidence interval:

Received April 2009; revised August 2009.
[1]Supported in part by NSF Grant DMS-09-06858.

(i) The accuracy: maintain the coverage probability function of the confidence interval at least $1 - \alpha$, that is,

(1) $\quad P_{(p_1,p_0)}(L(X,Y) \leq p_1 - p_0 \leq 1) \geq 1 - \alpha \quad \forall p_1, p_0 \in [0,1].$

Any interval satisfying (1) is called a one-sided $1 - \alpha$ confidence interval for $p_1 - p_0$. In general, there is no disagreement on criterion (i). When it is hard to implement (1), an approximate $1 - \alpha$ confidence interval is employed.

(ii) The precision: minimize the "size" (e.g., the length) of the confidence interval within a certain class of intervals.

Researchers do have different opinions on how to define an interval with the "minimum size." For two $1 - \alpha$ confidence intervals, $C_1(X,Y)$ and $C_2(X,Y)$, the most natural way to claim $C_1(X,Y)$ no worse than $C_2(X,Y)$ is that

(2) $\quad C_1(x,y)$ is a subset of $C_2(x,y)$ for all sample points $(x,y)$.

We call this the set inclusion criterion, proposed in Wang (2006), and an equivalent version was given in Bol'shev (1965). The superiority of $C_1$ over $C_2$ is easy to check because no expectation computation is involved. Under this criterion, for a specified class of $1 - \alpha$ intervals, we should search for the *smallest* $1 - \alpha$ confidence interval, which is equal to the intersection of all intervals in the class, provided that the intersection also belongs to the class.

For the case of the one-sided interval, the smallest interval in a class is equivalent to the shortest interval which has the shortest length on all sample points in that class. For the case of the two-sided interval, the smallest implies the shortest; however, the shortest does not necessarily imply the smallest. Also the smallest interval has a clear interpretation. Therefore, we use the terminology, "the smallest interval."

The existence of the smallest interval depends on the class of intervals from which we search for the smallest. In this paper, we propose two restrictions on the class:

(a) one-sided $1 - \alpha$ confidence interval;
(b) a certain monotonicity on the confidence limit $L(X,Y)$ for all sample points.

We will show the existence of the smallest interval and give its construction under these two restrictions. There were successful efforts for the case of a single proportion $p_1$ based on one observation $X$ where there exists a natural ordering on the lower limit $L(X)$: $L(x_1) \leq L(x_2)$ if $x_1 \leq x_2$, and there is no nuisance parameter. The smallest one-sided $1 - \alpha$ confidence interval for $p_1$ was derived independently by Bol'shev (1965) and Wang (2006). However, when there exists nuisance parameter, the result on the smallest one-sided confidence interval is very limited. Bol'shev and Loginov (1966) partially generalized Bol'shev's method (1965) to the case with nuisance parameter(s). Their construction is based on some function of $X$, $Y$ and $p_1 - p_0$ but not under condition (b). As we show later, for different orderings on $L(X,Y)$, the corresponding smallest intervals are different. So the interval following Bol'shev and Loginov's method for $p_1 - p_0$ cannot be the smallest.

The confidence interval was proposed by Laplace in 1814, and its definition only involves the coverage probability as shown in (1). However, the interval construction based on the coverage probability is not among the major methods currently used in practice. For example, Casella and Berger (1990), summarized five methods for the construction: the inversion of a family of tests; the pivotal quantities; a stochastic nonincreasing (or nondecreasing) distribution family with a single parameter, Bayesian intervals; invariant intervals, etc. The first is a general but indirect method because of inverting tests. During the inversion process, it is not easy to see how the interval is obtained. The other four need assumptions on the distribution under the study. For example, the second assumes the existence of pivotal quantities which is not true for binomial distributions. Being a major statistical inference procedure, the confidence interval deserves a direct method which is based on the analysis of coverage probability and needs mild or no assumptions on the distribution for its construction. However, it does need an assumption, restriction (b), on the sample space. The development of such a method is one goal of our paper. More importantly, this method can generate the smallest interval in many classes of intervals in the presence of nuisance parameter(s).

The rest of the paper is organized as follows. In Section 2, we specify appropriate classes of intervals, and in each class the smallest one-sided $1 - \alpha$ confidence interval for $p_1 - p_0$ is constructed. In Section 3, we carefully identify a special class of interval and then derive the corresponding smallest interval. An example is given to illustrate the procedure. In Section 4, we apply the interval in Section 3 to detect the minimum effective dose (MED), an important issue in clinical trials. A step-down test procedure is obtained with the familywise error rate controlled at level $\alpha$. In Section 5, we generalize the construction on the smallest one-sided confidence interval to other discrete sample spaces, and a Poisson distribution example is discussed. Some closing remarks are given in Section 6. A confidence interval with 0 confidence level is of no interest, so we assume $0 \leq \alpha < 1$ throughout the paper.

**2. The smallest one-sided confidence interval.** Recall $X \sim Bin(n, p_1)$ and $Y \sim Bin(m, p_0)$. Let $\Delta = p_1 - p_0$ be the parameter of interest and $p_0$ be the nuisance parameter. Let

(3) $\quad S = \{z = (x,y) : 0 \leq x \leq n, 0 \leq y \leq m, x \text{ and } y \text{ are integers}\}$

be the sample space. We use $z$ and $(x,y)$ interchangeably. Also the parameter space $H = \{(p_1, p_0) : 0 \leq p_1, p_0 \leq 1\}$ can be rewritten as

(4) $$H = \{(\Delta, p_0) : p_0 \in D(\Delta) \text{ for each } \Delta \in [-1,1]\},$$

where

(5) $$D(\Delta) = \begin{cases} [0, 1-\Delta], & \text{if } \Delta \in [0,1], \\ [-\Delta, 1], & \text{if } \Delta \in [-1,0). \end{cases}$$

The interval class that contains the smallest one-sided interval for $\Delta$ is given below.

DEFINITION 1. For any given ordered partition $\{C_j\}_{j=1}^{k_0}$ of $S$, define a class of one-sided $1-\alpha$ confidence intervals for $\Delta$;

$$\mathcal{B} = \{[L(Z), 1] : L(z) \text{ is constant on } C_j$$
$$\text{and } L(z) \geq L(z') \text{ for } z \in C_j, z' \in C_{j+1}, \forall j\}.$$

Since $L(z)$ is a constant on $C_j$, we define $L(C_j) = L(z)$ for any $z \in C_j$.

REMARK 1. Any given ordered partition of $S$ defines an ordering on the lower confidence limit. So we say: search for the smallest one-sided $1-\alpha$ confidence interval that is monotone with respect to the ordered partition $\{C_j\}_{j=1}^{k_0}$, or simply search for the smallest interval under the ordered partition. On the other hand, an ordered partition can be obtained by any given function $L(Z)$ as follows. Let $\{l_j\}_{j=1}^{k_0}$ be a sequence of strictly decreasing numbers in $j$ that contains all possible values of $L(Z)$. Then define $C_j = \{z \in S : L(z) = l_j\}$ for $1 \leq j \leq k_0$.

DEFINITION 2. A confidence interval $[L_S(Z), 1]$ in $\mathcal{B}$ is the smallest if it is a subset of any intervals in $\mathcal{B}$; that is, for any $[L(Z), 1]$ in $\mathcal{B}$, $L(z) \leq L_S(z), \forall z \in S$.

Definition 2 is adopted from Wang (2006). The smallest interval, if it exists, is the best in the strongest sense, and automatically minimizes the false coverage probability and the expected length among all intervals in $\mathcal{B}$. Next, we prove the existence and provide the construction of the smallest interval in $\mathcal{B}$.

THEOREM 1. *Assume $\alpha \in [0,1)$. For a given ordered partition $\{C_j\}_{j=1}^{k_0}$ of $S$ and any $z \in C_j$, let*

(6) $$f_j(\Delta) = \min_{p_0 \in D(\Delta)} P(S_j^c)$$
$$= \min_{p_0 \in D(\Delta)} \sum_{(x,y) \in S_j^c} p_X(x; p_0 + \Delta, n) p_Y(y; p_0, m),$$

*where* $S_j = \bigcup_{i=1}^{j} C_i$, *and let*

(7) $$G_z = \{\Delta \in [-1,1] : f_j(\Delta') \geq 1 - \alpha, \forall \Delta' < \Delta\}.$$

*Define*

(8) $$L_S(z) = \begin{cases} \sup G_z, & \text{if } G_z \neq \varnothing, \\ -1, & \text{otherwise.} \end{cases}$$

*Then:*

(1) $[L_S(Z), 1] \in \mathcal{B}$;
(2) $[L_S(Z), 1]$ *is the smallest in* $\mathcal{B}$.

REMARK 2. As pointed out in Corollary 1 in Wang (2006), the lower limit $L(x, y)$ is "the smallest $\theta$ ($\theta = \Delta$ in our case here) where $p(x; \theta)$ ($= p_X(x; n, \Delta + p_0) p_Y(y; m, p_0)$ in our case here) is used to compute the coverage probability." Therefore, in order to obtain the *largest* $L(x, y)$, we should exclude the term $p_X(x; n, \Delta + p_0) p_Y(y; m, p_0)$ from the coverage probability when $\Delta$ is as large as possible while keeping the coverage probability at least $1 - \alpha$. This is achieved by implementing (6)–(8) where (6) provides the minimal (respect to $p_0$) coverage probability (as a function of $\Delta$) for the desired interval up to $L_S(z)$, (7) assures its coverage probability no smaller than $1 - \alpha$ and (8) guarantees that $L_S(z)$ is the largest.

REMARK 3. Note that $f_j$, $G_z$ and $L_S$, defined in (6)–(8), depend on $z$ through $C_j$. Let $l_j = L_S(z)$ for $z \in C_j$.

PROOF OF THEOREM 1. To prove (1), first, it is clear that $L_S(z)$ is a constant on each $C_j$ due to Remark 3. Secondly, suppose $L_S(z) < L_S(z')$ for $z \in C_j$ and $z' \in C_{j+1}$ for some $j$. Then $L_S(z') > -1$. Notice

$$f_j(\Delta) \geq f_{j+1}(\Delta) \geq 1 - \alpha \qquad \forall \Delta < L_S(z').$$

So $L_S(z') \in G_z$ due to (7) and $L_S(z') \leq L_S(z)$ due to (8). A contradiction is constructed. So $L_S(z) \geq L_S(z')$ for $z \in C_j$ and $z' \in C_{j+1}$ for all $j$. Thirdly, consider the coverage probability function of $[L_S(Z), 1]$: $h_S(\Delta, p_0) = P(L_S(Z) \leq \Delta)$. Let

$$h(\Delta) = \min_{p_0 \in D(\Delta)} h_S(\Delta, p_0).$$

We need to show $h(\Delta) \geq 1 - \alpha$ on $[-1, 1]$. Note that $[-1, 1]$ is partitioned as $[l_1, 1] \cup (\bigcup_{j=1}^{k_0-1} [l_{j+1}, l_j))$ where $l_j$ is given in Remark 3. For $\Delta \in [l_1, 1]$, notice $L_S(z) \leq \Delta$ for all $z \in S$. Then

$$h(\Delta) = \min_{p_0 \in D(\Delta)} P(S) = 1 \geq 1 - \alpha.$$

For $\Delta \in [l_{j+1}, l_j)$ for any $1 \leq j \leq k_0 - 1$, notice $L_S(z) \leq \Delta$ for any $z \in S_j^c$. Then

$$h(\Delta) \geq \min_{p_0 \in D(\Delta)} P(S_j^c) = f_j(\Delta) \geq 1 - \alpha$$

by (7). Thus $[L_S(Z), 1] \in \mathcal{B}$.

To prove (2), suppose $[L_S(Z), 1]$ is not the smallest. Then there exists an interval $[L^*(Z), 1]$ in $\mathcal{B}$ and a point $z^* \in C_j$ for some $j$ with $L_S(z^*) < L^*(z^*)$. Let $l_i^* = L^*(z)$ for $z \in C_i$ for $i = 1, \ldots, k_0$. Then $l_j < l_j^*$. Let $h^*(\Delta, p_0)$ be the coverage probability of the interval $[L^*(Z), 1]$. For any $\Delta \in I = [l_j, l_j^*)$ (not an empty interval), we have

(9) $\quad 1 - \alpha \leq h^*(\Delta, p_0) = P(L^*(Z) \leq \Delta) \leq P(S_{j+1}^c) \leq P(S_j^c).$

The second inequality holds because $\{z : L^*(z) \leq \Delta\}$ is contained in $S_{j+1}^c$ when $\Delta \in I$. Therefore, $f_j(\Delta) = \min_{p_0 \in D(\Delta)} P(S_j^c) \geq 1 - \alpha$ on interval $I$ which contradicts (8). □

PROPOSITION 1. *For any one-sided $1 - \alpha$ confident limit $L(Z)$ with $0 \leq \alpha < 1$,*

(10) $$\min_{z \in S} L(z) = -1.$$

PROOF. Suppose $c = \min_{z \in S} L(z) > -1$. Pick a point $(\Delta_0, p_{00}) = ((-1 + c)/2, 1)$ in the parameter space. Note $\Delta_0 < L(z)$ for any $z \in S$; then

$$P_{(\Delta_0, p_{00})}(L(Z) \leq \Delta_0) = 0 < 1 - \alpha$$

which contradicts the fact that $[L(Z), 1]$ is of level $1 - \alpha$. □

EXAMPLE 1. Consider the case of $n = 4, m = 1$ and a predetermined ordered partition of $S$ given by the well-known $z$-test statistic,

$$Z(x, y) = \frac{\hat{p}_1 - \hat{p}_0}{\sqrt{\hat{p}_1(1 - \hat{p}_1)/n + \hat{p}_0(1 - \hat{p}_0)/m}},$$

following Remark 1 where $\hat{p}_1 = x/n$ and $\hat{p}_0 = y/m$, and $0/0 \stackrel{\text{def}}{=} 0$, $+/0 \stackrel{\text{def}}{=} \infty$ and $-/0 \stackrel{\text{def}}{=} -\infty$. Then this ordered partition $\{C_j^{ZT}\}$ and its associated smallest 95% confidence interval $[L_{ZT}(Z), 1]$ are reported in Table 1.

For the purpose of illustration, we determine $L_{ZT}(3, 0)$ here. Consider

$$f_2(\Delta) = \min_{p_0 \in D(\Delta)} (1 - p_X(4; 4, \Delta + p_0)) p_Y(0; 1, p_0)$$
$$- p_X(3; 4, \Delta + p_0) p_Y(0; 1, p_0)).$$

Since $f_2(-1) = 1$, $L_{ZT}(3,0)$ is equal to $-0.345$, the smallest solution of $f_2(\Delta) = 1 - \alpha$ with $\alpha = 0.05$. This can be done numerically by calculating $f_2(\Delta)$ at each $\Delta$ in the order of $\Delta = -1, -0.999, -0.998, \ldots$ with an increment of 0.001, for example, until $f_2(\Delta)$ is always greater than $1 - \alpha$. Therefore, the last value of $\Delta = -0.345$ is the smallest solution.

It is well known that a $1 - \alpha$ confidence interval can generate a family of level-$\alpha$ tests and vice versa. Interval $[L_S(Z), 1]$ in Theorem 1 can be used for this purpose for the following family of hypotheses:

$$H_0(\delta): \Delta \leq \delta \quad \text{vs.} \quad H_A(\delta): \Delta > \delta, \tag{11}$$

where $\delta \in [-1, 1]$. For any given $\delta$, the rejection region,

$$R_S(\delta) = \{z \in S : L_S(z) > \delta\}, \tag{12}$$

defines a level-$\alpha$ test for (11). For the ordered partition $\{C_j\}_{j=1}^{k_0}$ of $S$, let

$$j(\delta) = \max\{j : L_S(C_j) > \delta\};$$

or $j(\delta) = 0$, if $L_S(C_1) \leq \delta$. Then

$$R_S(\delta) = \bigcup_{j=1}^{j(\delta)} C_j$$

due to $L_S(C_{j+1}) \leq L_S(C_j)$ and the definition of $j(\delta)$. On the other hand, $[L_S(Z), 1]$ can also be obtained by inverting tests as follows. For an ordered partition $\{C_j\}_{j=1}^{k_0}$, consider a level-$\alpha$ rejection region of form $\bigcup_{j=1}^{s(\delta)} C_j$ for

TABLE 1
*Different ordered partitions and their associated smallest 95% intervals when $n = 4$ and $m = 1$*

| $j$ | $C_j^{ZT}$ | $Z$ | $L_{ZT}(z)$ | $j$ | $C_j^Z$ | $L_A(C_j^Z)$ | $L_Z(z)$ | $j$ | $C_j^I$ | $L_I(z)$ |
|---|---|---|---|---|---|---|---|---|---|---|
| 1 | (4, 0) | $\infty$ | $-0.095$ | 1 | (4, 0) | 1 | $-0.095$ | 1 | (4, 0) | $-0.095$ |
| 2 | (3, 0) | 3.464 | $-0.345$ | 2 | (3, 0) | 0.394 | $-0.345$ | 2 | (3, 0) | $-0.345$ |
| 3 | (2, 0) | 2.000 | $-0.562$ | 3 | (2, 0) | 0.089 | $-0.562$ | 3 | (2, 0) | $-0.562$ |
| 4 | (1, 0) | 1.155 | $-0.756$ | 4 | (4, 1) | 0 | $-0.950$ | 4 | (1, 0) | $-0.756$ |
| 5 | (4, 1) | 0 | $-0.950$ | — | (0, 0) | 0 | — | 5 | (4, 1) | $-0.757$ |
|  | (0, 0) | 0 | — | 5 | (1,0) | $-0.106$ | $-0.950$ | 6 | (3, 1) | $-0.770$ |
| 6 | (3, 1) | $-1.155$ | $-0.950$ | 6 | (3, 1) | $-0.606$ | $-0.950$ | 7 | (2, 1) | $-0.902$ |
| 7 | (2, 1) | $-2.000$ | $-0.950$ | 7 | (2, 1) | $-0.911$ | $-0.950$ | 8 | (0, 0) | $-0.950$ |
| 8 | (1, 1) | $-3.464$ | $-0.987$ | 8 | (0, 1) | $-1$ | $-1$ | 9 | (1, 1) | $-0.987$ |
| 9 | (0, 1) | $-\infty$ | $-1$ | 5 | (1, 1) | $-1.106$ | $-1$ | 10 | (0, 1) | $-1$ |

some nonnegative integer $s(\delta)$ for hypothesis $H_0(\delta)$ given in (11) with a fixed $\delta$ where

$$s(\delta) = \max\left\{n \leq k_0 : \sup_{(\Delta,p_0) \in H_0(\delta)} P\left(\bigcup_{j=1}^{n} C_j\right) \leq \alpha\right\}. \tag{13}$$

So $(\bigcup_{j=1}^{s(\delta)} C_j)^c$, the complement of $\bigcup_{j=1}^{s(\delta)} C_j$, is the acceptance region. For a sample point $Z = z$, let

$$C_T(z) = \left\{\delta \in [-1,1] : z \in \left(\bigcup_{j=1}^{s(\delta)} C_j\right)^c\right\}. \tag{14}$$

Then $[L_S(Z), 1]$ equals $C_T(Z)$ as shown in Theorem 2 below. However, we lose the intuition given in (6)–(8) during this inversion process.

THEOREM 2. *$C_T(Z)$ belongs to the interval class $\mathcal{B}$ given in Definition 1, and*

$$[L_S(Z), 1] = C_T(Z). \tag{15}$$

PROOF. First, for any sample point $z$, if $\Delta \in C_T(z)$ and $\Delta' \in [\Delta, 1]$, then $s(\Delta') \leq s(\Delta)$ following (13) and (14). Therefore, $\Delta' \in C_T(z)$, and $C_T(z)$ is a confidence interval for $\Delta$. Second, the coverage probability of $C_T(Z)$

$$P(\Delta \in C_T(Z)) = P(z : \Delta \in C_T(z)) = P\left(z \in \left(\bigcup_{j=1}^{s(\Delta)} C_j\right)^c\right) \geq 1 - \alpha$$

following (13) for any given $(\Delta, p_0)$. So $C_T(Z)$ is of level $1 - \alpha$. Third, let $C_T(z) = [L_T(z), 1]$. (i) It is clear that $L_T(z)$ is constant on each $C_i$; (ii) Pick $z_1 \in C_i$, any $\delta_1 \in C_T(z_1)$ and $z_2 \in C_{i+1}$. Since $z_1 \in (\bigcup_{j=1}^{s(\delta_1)} C_j)^c$, we have $s(\delta_1) < i$. Thus $s(\delta_1) < i + 1$ and $z_2 \in (\bigcup_{j=1}^{s(\delta_1)} C_j)^c$, and we conclude that $\delta_1 \in C_T(z_2)$. Therefore, $L_T(C_i)$ is nonincreasing in $i$ for $1 \leq i \leq k_0$. So $C_T(Z)$ belongs to the interval class $\mathcal{B}$. This implies $[L_S(Z), 1] \subset C_T(Z)$ following Theorem 1.

Now we only need to prove

$$[L_S(Z), 1] \supset C_T(Z)(= [L_T(Z), 1]). \tag{16}$$

Without loss of generality, assume $L_S(C_j)$ is strictly decreasing in $j$. Otherwise we redefine a new ordered partition $\{C'_j\}$, by merging those $C_j$s on which $L_S(C_j)$ does not change so that $L_S(C'_j)$ is strictly decreasing. Suppose (16) is not true. Then let $j_0$ be the smallest positive integer so that $L_T(C_{j_0}) <$ $L_S(C_{j_0})$. Pick $\delta_0 \in (\max\{L_T(C_{j_0}), L_S(C_{j_0+1})\}, L_S(C_{j_0}))$. Note $\bigcup_{j=1}^{s(\delta_0)} C_j = \{z : L_T(z) \leq \delta_0\}^c = \bigcup_{j=1}^{j_0-1} C_j$. So

(17) $$s(\delta_0) = j_0 - 1.$$

On the other hand,

$$1 - \alpha \leq \inf_{(\Delta, p_0) \in H_0(\delta_0)} P(L_S(Z) \leq \Delta) \leq \inf_{(\Delta, p_0) \in H_0(\delta_0)} P(L_S(Z) \leq \delta_0)$$

$$= 1 - \sup_{(\Delta, p_0) \in H_0(\delta_0)} P\left(\bigcup_{j=1}^{j_0} C_j\right).$$

Due to (13), $s(\delta_0) \geq j_0$, which contradicts (17). Then, (16) is true, as well as (15). □

REMARK 4. When applying Theorem 1, our intention is not to generate the optimal interval among all possible orderings, but to improve or modify a given interval $[L(Z), 1]$ which has a level $1 - \alpha$ or approximately $1 - \alpha$, to be the smallest $1 - \alpha$ interval. To achieve this, one forms an ordered partition for $S$ following Remark 1 using the given function $L(Z)$, then derives the smallest interval $[L_S(Z), 1]$ following Theorem 1.

EXAMPLE 1 (Continued). The most commonly used one-sided interval for $\Delta$ in practice is the following $z$-interval:

(18) $$[L_A(z), 1] \stackrel{\text{def}}{=} [\hat{p}_1 - \hat{p}_0 - z_\alpha \sqrt{\hat{p}_1(1-\hat{p}_1)/n + \hat{p}_0(1-\hat{p}_0)/m}, 1],$$

where $z_\alpha$ is the upper $\alpha$th percentile of the standard normal distribution. Its coverage probability can be much less than the nominal level $1 - \alpha$. We follow Remark 4 to modify this interval by generating an ordered partition of $S$, denoted by $\{C_j^Z\}$. Then the smallest $1 - \alpha$ interval, denoted by $[L_Z(Z), 1]$, based on this partition is derived following Theorem 1 and is reported in Table 1 for the case in Example 1. Note $C_1^Z = (1, 1)$ instead of $(0,1)$ which is intuitively incorrect. Therefore, interval $[L_Z(Z), 1]$ is not recommended.

Is it possible to improve the smallest interval from a given ordered partition? Yes, especially when there exists a finer partition than the given one as stated below. See such an example in Table 1 where $L_I(Z) \geq L_{ZT}(Z) \geq L_Z(Z)$. Each $C_j^I$, given in the second-to-last column of Table 1, contains a single sample point which implies that $\{C_j^I\}$ is a finest partition of $S$.

PROPOSITION 2. *For two ordered partitions* $\mathcal{P} = \{C_j\}_{j=1}^{k_0}$ *and* $\mathcal{P}^* = \{C_j^*\}_{j=1}^{k_0^*}$ *of the sample space* $S$, *suppose each* $C_j$ *is a subset of* $C_{i(j)}^*$ *for some* $i(j)$ where $i(j)$ is a nondecreasing function in $j$ (i.e., $\mathcal{P}$ is a finer partition than $\mathcal{P}^*$). Let $[L_S(Z), 1]$ and $[L_S^*(Z), 1]$ be the smallest $1 - \alpha$ confidence intervals under ordered partitions $\mathcal{P}$ and $\mathcal{P}^*$, respectively. Then

(19) $$L_S^*(Z) \leq L_S(Z) \qquad \forall z \in S.$$

PROOF. The proof is trivial if one notices that any ordering on $L(Z)$ by $\mathcal{P}^*$ is also an ordering by $\mathcal{P}$. Then the claim follows Theorem 1. □

**3. An ordering on the confidence limits.** Which ordered partition of $S$ provides an interval that cannot be uniformly improved? Roughly speaking, by Theorem 1, we prefer an ordering on $L(z)$ $(= L(x, y))$

(20) that yields a large smallest solution of $f_j(\Delta) = 1 - \alpha$ for all $j$s.

Due to Proposition 2, each set in the partition would contain only one point. Because of the specialty of binomial distributions, $L(x, y)$ should satisfy:

(1) $L(x_1, y) \leq L(x_2, y)$ for $x_1 \leq x_2$; and
(2) $L(x, y_1) \geq L(x, y_2)$ for $y_1 \leq y_2$.

Let $B_B$ denote the class of all one-sided $1 - \alpha$ intervals for $\Delta$ satisfying (1) and (2). We will search for optimal intervals, perhaps admissible ones, from $B_B$ in this section.

It is clear that $L(n, 0)$ must be the largest among all $L(x, y)$s and the second largest $L(x, y)$ should be achieved at either $(n - 1, 0)$ or $(n, 1)$ or both (if $n = m$). We, by induction, construct an ordered partition of $S$, denoted by $\{C_j^B\}_{j=1}^{k_0^B}$, that satisfies (1) and (2) and starts at point $(n, 0)$ as follows:

Step 1: Let $C_1^B = \{(n, 0)\}$, $m_1 = 1$ and $m_0 = 0$ because $L(n, 0)$ is the largest among all $L(x, y)$s. So $C_1^B = \{(x_i, y_i)\}_{i=m_0+1}^{m_1}$ where $(x_1, y_1) = (n, 0)$.

Step 2: Suppose, by induction, $\{C_j^B\}_{j=1}^{k}$ are available for some positive integer $k$ where

$$C_j^B = \{(x_i, y_i)\}_{i=m_{j-1}+1}^{m_j}$$

for some nonnegative integers $m_0, m_1, \ldots, m_k$ satisfying:

(I) $L(x, y)$ is constant on $C_j^B$;
(II) $L(x_{m_{j-1}}, y_{m_{j-1}}) \geq L(x_{m_j}, y_{m_j})$ for each $j \leq k$.

Now we determine $C_{k+1}^B$. Let $S_k = \bigcup_{j=1}^{k} C_j^B$ and let $N_k$ be the "neighbor" set of $S_k$, that is,

$$N_k = \{(x, y) \in S : (x, y) \notin S_k; (x+1, y) \in S_k \text{ or } (x, y-1) \in S_k\}.$$

Due to (1) and (2), some points in $N_k$ are disqualified to be in $C_{k+1}^B$. To exclude these points, let $NC_k$ be the "candidate" set within $N_k$ satisfying

(21) $\qquad NC_k = \{(x,y) \in N_k : (x+1, y) \notin N_k \text{ and } (x, y-1) \notin N_k\}.$

Therefore, $C_{k+1}^B$ must be a subset of $NC_k$, and a point selected from $NC_k$ automatically guarantees (1) and (2). For each point $z_0 = (x_0, y_0)$ in $NC_k$, consider

$$f_{z_0}(\Delta) = \min_{p_0 \in D(\Delta)} P((\{z_0\} \cup S_k)^c)$$

$$= \min_{p_0 \in D(\Delta)} \sum_{z \in (\{z_0\} \cup S_k)^c} p_X(x; p_0 + \Delta, n) p_Y(y; p_0, m).$$

Let

(22) $\qquad E_{z_0} = \{\Delta \in [-1, 1] : f_{z_0}(\Delta') \geq 1 - \alpha, \forall \Delta' < \Delta\}$

and

(23) $\qquad L_o(z_0) = \begin{cases} \sup E_{z_0}, & \text{if } E_{z_0} \neq \varnothing, \\ -1, & \text{otherwise.} \end{cases}$

Define

(24) $\qquad C_{k+1}^B = \left\{ z \in NC_k : L_o(z) = \max_{z_0 \in NC_k} L_o(z_0) \right\}$ and

(25) $\qquad m_{k+1} = m_k + \text{the number of elements in } C_{k+1}^B.$

Note that $C_{k+1}^B$ may contain more than one point especially when $n = m$. By induction, an ordered partition $\{C_j^B = \{z_i\}_{i=m_{j-1}+1}^{m_j}\}_{j=1}^{k_0^B}$ for $S$ with some positive integer $k_0^B$ is constructed. Therefore, the smallest one-sided $1 - \alpha$ confidence interval under this ordered partition, denoted by $[L_S(Z), 1]$, is constructed for estimating $\Delta$ following Theorem 1.

REMARK 5. $E_{z_0}$ and $G_z$ ($f_{z_0}$ and $f_j$, $L_o$ and $L_S$) are defined in a similar way. From (24), the ordered partition $\{C_j^B\}_{j=1}^{k_0^B}$ tends to yield a large $L_o(z)$, which results in a short interval compared with other partitions. More precisely, $L_S(C_j^B)$ equals the largest possible value provided that $L_S(C_1^B), \ldots, L_S(C_{j-1}^B)$ are determined. However, different from $(G_z, f_j, L_S)$ that depends on $z$ through $C_j^B$, $(E_{z_0}, f_{z_0}, L_o)$ depends on each individual $z_0$. If $C_j^B$ always contains a single point for any $j$, then, for $z \in C_j^B$, $L_S(z)$ equals the largest of $L_o(z_0)$s in the previous step, and we have the following result.

PROPOSITION 3. *For ordered partition $\{C_j^B\}$ and interval $[L_S(Z), 1]$ constructed in steps 1 and 2, if each $C_j^B$ contains only one sample point (i.e., $m_j = m_{j-1} + 1$ for all $j$s), then $[L_S(Z), 1]$ is admissible in $B_B$. That is, for an interval $[L(Z), 1] \in B_B$, if $L_S(z) \leq L(z)$ for any $z \in S$, then $L_S(z) = L(z)$ for any $z \in S$.*

PROOF. Suppose the claim is not true. Note each $C_j^B$ containing only one point, and let $j_0$ be the smallest positive integer so that $L_S(C_{j_0}^B) < L(C_{j_0}^B)$. Let $\{C_j\}$ be the associated ordered partition for $[L(Z), 1]$. Then $C_j^B = C_j$ for any $j < j_0$. Therefore, $C_{j_0}$ is a subset of $NC_{j_0-1}$ given in (21) due to conditions (1) and (2). Noting (24), we conclude $L_S(C_{j_0}^B) \geq L(C_{j_0}^B)$. A contradiction is constructed. □

Conditions (1) and (2) were first proposed in Barnard (1947) and called the "C" condition. He constructed an optimal rejection region for a hypothesis testing problem,

$$H_0(0) : \Delta \leq 0 \quad \text{vs.} \quad H_A(0) : \Delta > 0,$$

a special case of (11) for $\delta = 0$, using a special ordering on $S$. This ordering satisfies conditions (1) and (2), and the corresponding ordered partition $\{C_j\}_{j=1}^{k_0}$ is generated by induction starting at $C_1 = (n, 0)$. Also, for given $C_1, \ldots, C_{j_0-1}$ with a positive integer $j_0$ ($\leq k_0$), $C_{j_0}$ is chosen so that $\sup_{p_0 \in D(0)} P(\bigcup_{j=1}^{j_0} C_j)$ is minimized. So Barnard's ordering is similar to ours, except that he focused on $\Delta = 0$, but we deal with all $\Delta \in [-1, 1]$. Pointed out by Martin Andres and Silva Mato (1994), Barnard's test is the (overall) most powerful existing test for comparing two independent proportions.

EXAMPLE 1 (Continued). Now construct the smallest 95% confidence interval with partition $\{C_j^B\}_{j=1}^{k_0^B}$. First $C_1^B = \{(4, 0)\}$ following step 1, and $L_S(4, 0) = -0.095$ by solving

$$f_1(\Delta) = \min_{\{p_0 \in D(\Delta)\}} (1 - p_X(4; 4, \Delta + p_0) p_Y(0; 1, p_0)) = 0.95$$

because $f_1$ now is nonincreasing in $\Delta$. In step 2, $N_1$, the neighbor set of $S_1 (= C_1^B)$, is equal to $\{(3, 0), (4, 1)\}$, and $NC_1 = N_1$. Following (23),

$$L_o(3, 0) = -0.345, \qquad L_o(4, 1) = -0.527.$$

Thus $C_2^B = \{(3, 0)\}$ by (24). In step 3, three sets are needed, $S_2$, $N_2$ and $NC_2$, and are given below in the sample space $S$. Note here $NC_2 \neq N_2$.

|   |   |   | $x$ |   |   |
|---|---|---|---|---|---|
| $y$ | 0 | 1 | 2 | 3 | 4 |
| 1 | – | – | – | $N_2$ | $N_2$, $NC_2$ |
| 0 | – | – | $N_2$, $NC_2$ | $S_2$ | $S_2$ |

Again, for each point in $NC_2$, we have $L_o(2,0) = -0.561, L_o(4,1) = -0.527$, following (23). Then $C_3^B = \{(4,1)\}$ by (24). The rest of the interval construction is given in Table 2. Following Remark 5, since each $C_j^B$ contains a single point, $L_S(z)$ on $C_j^B$ is equal to the largest $L_o(z_0)$ in the previous step and is reported in the last column of Table 2, and the construction is complete at the 10th ($= k_0^B$) step. This interval is admissible in $B_B$ due to Proposition 3. However, if compared with interval $[L_I(Z), 1]$ in Table 1, neither uniformly dominates the other.

**4. Identifying the minimum effective dose.** Suppose we have a sequence of independent binomial random variables $X_i \sim Bin(n_i, p_i)$ for $i = 1, \ldots, k$ and $Y \sim Bin(m, p_0)$. The goal here is to identify the smallest positive integer $i_0$ so that $p_i > p_0 + \delta$ for any $i \in [i_0, k]$, where $\delta$ is some predetermined nonnegative number. Each $p_i$ is the proportion of patients who show improvement using a drug at dose level $i$. A large $i$ associates with a large dose level, and $p_0$ is the proportion for the control group. Then $i_0$ is called the minimum effective dose (MED). Finding the MED is important since high doses often turn out to have undesirable side effects.

TABLE 2
The details of the construction of partition $\{C_j^B\}_{j=1}^{k_0^B}$ when $n = 4$ and $m = 1$

| $j$ | $C_j^B$ | $N_j$ | $NC_j$, $L_o(z_0)$ | $L_S(C_j^B)$ |
|---|---|---|---|---|
| 1 | (4,0) | (3,0), (4,1) | (3,0), (4,1)<br>**−0.345**, −0.527 | −0.095 |
| 2 | (3,0) | (2,0), (3,1), (4,1) | (2,0), (4,1)<br>−0.561, **−0.527** | **−0.345** |
| 3 | (4,1) | (2,0), (3,1) | (2,0), (3,1)<br>**−0.578**, −0.752 | **−0.527** |
| 4 | (2,0) | (1,0), (2,1), (3,1) | (1,0), (3,1)<br>−0.757, **−0.752** | **−0.578** |
| 5 | (3,1) | (1,0), (2,1) | (1,0), (2,1)<br>**−0.770**, −0.902 | **−0.752** |
| 6 | (1,0) | (0,0), (1,1), (2,1) | (0,0), (2,1)<br>−0.950, **−0.902** | **−0.770** |
| 7 | (2,1) | (0,0), (1,1) | (0,0), (1,1)<br>**−0.950**, −0.987 | **−0.902** |
| 8 | (0,0) | (0,1), (1,1) | (1,1)<br>**−0.987** | **−0.950** |
| 9 | (1,1) | (0,1) | (0,1)<br>**−1** | **−0.987** |
| 10 | (0,1) |  |  | **−1** |

Typically, the MED is to be found when $X_i$ follows a normal distribution with the comparison in proportions replaced by that in means. Thus, the assumption of normality is an issue to be addressed. See, for example, Tamhane, Hochberg and Dunnett (1996), Hsu and Berger (1999), Bretz, Pinheiro and Branson (2005) and Wang and Peng (2008) for results under this setting. Now we search for the MED with a binary response without such concern on the distribution; see Tamhane and Dunnett (1999).

A sequence of hypotheses can be formulated to detect the MED as follows:

$$
\begin{aligned}
&H_{0i}: \min_{j \geq i}\{p_j - p_0\} \leq \delta \quad \text{vs.} \\
&H_{Ai}: \min_{j \geq i}\{p_j - p_0\} > \delta, \qquad \text{for } i = 1, \ldots, k,
\end{aligned}
\tag{26}
$$

which is similar to the one in Hsu and Berger (1999), page 471. It is clear that the MED equals the smallest $i$ for which $H_{Ai}$ is true. Also $H_{0i}$ is decreasing in $i$, thus $\mathcal{C} = \{H_{0i} : i = 1, \ldots, k\}$ is closed under the operation of intersection. Suppose a level $\alpha$ nondecreasing (in $i$) rejection region $R_i$ for $H_{0i}$ is constructed. Then, for the multiple test problem for testing all null hypotheses in $\mathcal{C}$, define a multiple test procedure: assert $H_{Ai}$ if $R_i$ occurs. This procedure controls the familywise error rate at level $\alpha$ following the closed test procedure by Marcus, Peritz and Gabriel (1976).

Now we apply the interval derived in Section 3 to obtain a level $\alpha$ test for $H_{0i}$. Let $L_{S,i}(X_i, Y)$ be the smallest one-sided $1-\alpha$ confidence interval for $p_i - p_0$ obtained in Section 3 before Remark 5. Define a rejection region for $H_{0i}$:

$$
R_i = \left\{(x_1, \ldots, x_k, y) : \min_{i \leq j \leq k}\{L_{S,j}(x_j, y)\} > \delta\right\}.
\tag{27}
$$

THEOREM 3. *The rejection region $R_i$ is nondecreasing in $i$ and is of level $\alpha$ ($<1$) for $H_{0i}$. Therefore, the multiple test procedure, which asserts not $H_{0i}$ (i.e., asserts $H_{Ai}$) if $R_i$ occurs for any $H_{0i} \in \mathcal{C}$, controls the familywise error rate at level $\alpha$.*

PROOF. First, it is trivial that $R_i$ is nondecreasing in $i$. Secondly, for any $(p_1, \ldots, p_k, p_0) \in H_{0i}$, there exists an $i^* \in [i, k]$ satisfying $p_{i^*} - p_0 \leq \delta$. Then $P_{(p_{i^*}, p_0)}(L_{S,i^*}(X_{i^*}, Y) > \delta) \leq P_{(p_{i^*}, p_0)}(L_{S,i^*}(X_{i^*}, Y) > p_{i^*} - p_0)$ when $p_{i^*} - p_0 \leq \delta$. Since $P_{(p_{i^*}, p_0)}(L_{S,i^*}(X_{i^*}, Y) > p_{i^*} - p_0) \leq \alpha$, we have

$$P_{(p_1, \ldots, p_k, p_0)}(R_i) \leq P_{(p_{i^*}, p_0)}(L_{S,i^*}(X_{i^*}, Y) > \delta) \leq \alpha.$$

The rest of the theorem follows the closed test procedure by Marcus, Peritz and Gabriel (1976) because $R_i$ is nondecreasing and $H_{0i}$ is decreasing in $i$. □

REMARK 6. The multiple test procedure with rejection regions $\{R_i\}_{i=1}^k$ in (27) is equivalent to the following step-down test procedure.

*Step* 1. If $R_k$ does not occur, conclude that the MED does not exist and stop; otherwise go to the next step.

*Step* 2. If $R_{k-1}$ does not occur, conclude the $MED = k$ and stop; otherwise go the next step.

$\vdots$

*Step k*. If $R_1$ does not occur, conclude the $MED = 2$ and stop; otherwise conclude the $MED = 1$ and stop.

REMARK 7. The proposed multiple test procedure is valid without the assumption of $p_1 \leq p_2 \leq \cdots \leq p_k$.

**5. A generalization.** Suppose a random vector $\underline{X}$ is observed from a discrete sample space $S$ with either finite or countable sample points, that is, $S = \{\underline{x}_i\}_{i=a}^b$ where $-\infty \leq a < b \leq +\infty$. An ordered partition $\{C_j\}_{j=c}^d$ on $S$ is given for some $-\infty \leq c \leq d \leq +\infty$. The probability mass function of $\underline{X}$ is given by $p(\underline{x}; \underline{\theta})$ where $\underline{\theta}$ is the parameter vector belonging to a parameter space $\Theta$, a subset of $R^k$. Suppose $\underline{\theta} = (\theta, \underline{\eta})$ and

$$\Theta = \{\underline{\theta} : \underline{\eta} \in D(\theta) \text{ for each } \theta \in [A, B]\},$$

where $[A, B]$ is a given interval in $R^1$(A and B may be $\pm\infty$, and the interval is open when the corresponding ending is infinity), and $D(\theta)$ is a subset of $R^{k-1}$ depending on $\theta$. Now we are interested in searching for the smallest one-sided $1 - \alpha$ confidence interval of form $[L(\underline{X}), B]$ for $\theta$ under the ordered partition $\{C_j\}_{j=c}^d$, i.e., $L(\underline{x})$ is constant on each $C_j$ and $L(\underline{x}) \geq L(\underline{x}')$ for any $\underline{x} \in C_j$ and $\underline{x}' \in C_{j'}$ for any $j \leq j'$.

THEOREM 4. *Assume $\alpha \in [0, 1)$. For a given partition $\{C_j\}_{j=c}^d$ of $S$ and any $\underline{x} \in C_j$, let*

$$(28) \quad f_j(\theta) = \inf_{\underline{\eta} \in D(\theta)} P(S_j^c) = \inf_{\underline{\eta} \in D(\theta)} \sum_{\underline{z} \in S_j^c} p(\underline{z}; \theta, \underline{\eta}),$$

*where $S_j = \bigcup_{i=c}^j C_i$ and let*

$$(29) \quad GG_{\underline{x}} = \{\theta \in [A, B] : f_j(\theta') \geq 1 - \alpha, \forall \theta' < \theta\}.$$

*Define*

$$(30) \quad L_G(\underline{x}) = \begin{cases} \sup GG_{\underline{x}}, & \text{if } GG_{\underline{x}} \neq \varnothing, \\ A, & \text{otherwise.} \end{cases}$$

*Then $[L_G(\underline{X}), B]$ is the smallest one-sided $1 - \alpha$ confidence interval under partition $\{C_j\}_{j=c}^d$.*

The proof is similar to Theorem 1 and is omitted.

EXAMPLE 2. Suppose one is interested in the difference, $\Delta$, of two means, $\lambda_1$ and $\lambda_2$, of two independent Poisson random variables, $X$ and $Y$. The sample space $S$ and the parameter space are

$S = \{(x,y) : x \text{ and } y \text{ are nonnegative integers}\}$ and

$H = \{(\Delta, \lambda_2) : \lambda_2 \in [0, +\infty) \text{ if } \Delta \in [0, +\infty); \lambda_2 \in [-\Delta, +\infty) \text{ if } \Delta \in (-\infty, 0)\}$,

respectively. One-sided $1 - \alpha$ confidence intervals of form $[L(X,Y), +\infty)$ for $\Delta$ are of interest. Let $F_P(x; \lambda)$ be the cumulative distribution function of a Poisson distribution with mean $\lambda$.

Different from the binomial case, there exists no fixed sample point on which $L(x,y)$ is the largest or the smallest. Lack of a starting or an ending point, a construction of an ordered partition $\{C_j\}_{j=c}^d$ of $S$ by induction is difficult. Instead, we show how to improve a naive interval $[L_1(X,Y), \infty)$ given below for $\Delta$.

This interval is obtained combining two smallest one-sided intervals $[L(X), +\infty)$ for $\lambda_1$ and $[0, U(Y)]$ for $\lambda_2$, both of level $\sqrt{1-\alpha}$. Following Bol'shev (1965), $L(x)$ satisfies

(31) $\quad 1 - F_P(x-1; L(x)) = 1 - \sqrt{1-\alpha} \qquad$ for $x > 0$ and $L(0) = 0$,

and $U(y)$ satisfies

(32) $$F_P(y; u(y)) = 1 - \sqrt{1-\alpha}.$$

Then $[L_1(X,Y), \infty)$, where $L_1(X,Y) \stackrel{\text{def}}{=} L(X) - U(Y)$, is a one-sided confidence interval of level $1 - \alpha$ for $\Delta = \lambda_1 - \lambda_2$ because

$P(L_1(X,Y) \leq \Delta) \geq P(L(X) \leq \lambda_1, U(Y) \geq \lambda_2) \geq (\sqrt{1-\alpha})^2 = 1 - \alpha.$

It is clear that $L_1(x,y)$ satisfies (1) and (2) introduced at the beginning of Section 3 because $L(x)$ and $U(y)$ are both increasing functions. Now we improve $[L_1(X,Y), +\infty)$ by constructing $L_G(X,Y)$ following Theorem 4 and Remark 4.

To illustrate the procedure, suppose $(X,Y) = (4,2)$ is observed and $\alpha = 0.05$. Following (31) and (32), we obtain $L(4) = 1.094$ and $U(2) = 7.208$, respectively, and $L_1(4,2) = -6.114$. We need to determine $L_G(4,2)$ given in (30). Consider a subset $S_j$ of $S$ on which $L_1(x,y)$ is no smaller than $L_1(4,2)$, that is,

$S_j = \{(x,y) \in S : L_1(x,y) \geq L_1(4,2)\} = \{(x,y) \in S : x \geq g(y)\},$

where, for each $y$,

$g(y) \stackrel{\text{def}}{=} \min\{x \geq 0 : L_1(x,y) \geq L_1(4,2)\}.$

For example, $g(0) = g(1) = 0$, and $g(3) = 7$. Plug $S_j$ into (28) and solve $L_G(4,2) = -4.744$ following (29) and (30). $L_G(4,2)$ is much larger than $L_1(4,2)$, as is well expected.

**6. Discussion.** In this paper, we discuss how to derive the smallest confidence interval under an ordered partition for a parameter in the presence of nuisance parameter(s) when the sample space is discrete. The interval construction is based on a direct analysis on the coverage probability, and needs an ordering on the sample points and mild assumptions (e.g., a discrete sample space) on the underlying distribution. The set inclusion criterion is employed for searching for good intervals because it has a clear interpretation. Under this criterion, the smallest interval is the best in the strongest sense provided its existence. It is well known that the existence of the best interval depends on the class of intervals from which the best is searched. We successfully characterize such classes by (a) considering one-sided $1 - \alpha$ confidence intervals and (b) requiring an ordering on the random confidence limits. Bol'shev and Loginov (1966) did not construct the interval under (b), so their method typically does not generate the smallest interval, while ours does. Another application of the proposed method is to identify admissible confidence intervals more efficiently. As an example, consider the case in Section 2. Let $[L_{S,C}(X,Y), 1]$ be the smallest $1 - \alpha$ confidence interval for $p_1 - p_0$ corresponding to a partition $\{C_j\}_{j=1}^{k_0}$ of $S$. Then the class of $1 - \alpha$ confidence intervals,

$$\mathcal{D} = \{[L_{S,C}(X,Y), 1] : \forall \text{ partition of } S\},$$

is complete since for any $1 - \alpha$ confidence interval $[L(X,Y), 1]$ there exists an interval $[L_{S,C}(X,Y), 1]$ in $\mathcal{D}$ so that $[L_{S,C}(x,y), 1]$ is always a subset of $[L(x,y), 1]$ for any $(x,y) \in S$. Although class $\mathcal{D}$ is not minimal, it contains finitely many elements. One only needs to search for optimal intervals from $\mathcal{D}$ because it is complete. Furthermore, one can apply Proposition 2 and conditions (1) and (2) in Section 3 to search for optimal intervals from a much smaller subset, also complete in $B_B$, of $\mathcal{D}$.

**Acknowledgments.** The author is grateful to three anonymous referees and an Associate Editor for their constructive suggestions. The author also thanks Professor Thaddeus Tarpey for proof reading the manuscript.

Department of Mathematics and Statistics
Wright State University
Dayton, Ohio 45435
USA
E-mail: weizhen.wang@wright.edu